\documentclass[11pt,reqno]{amsart}
\usepackage{fullpage}

\usepackage{amsmath,amssymb,amsthm,amsfonts,hhline,color}
\usepackage{etoolbox, textcomp}

\usepackage{hyperref}
\definecolor{my-linkcolor}{rgb}{0.75,0,0}
\definecolor{my-citecolor}{rgb}{0.1,0.57,0}
\definecolor{my-urlcolor}{rgb}{0,0,0.75}
\hypersetup{
	colorlinks, 
	linkcolor={my-linkcolor},
	citecolor={my-citecolor}, 
	urlcolor={my-urlcolor}
}

\def\Z{{\mathbb Z}}

\def\QQ{{\mathbb Q}}

\DeclareMathOperator{\even}{even}
\DeclareMathOperator{\odd}{odd}
\DeclareMathOperator{\sdim	}{sdim}
\DeclareMathOperator{\osp}{\mathfrak{osp}}
\DeclareMathOperator{\ch}{ch} 
\DeclareMathOperator{\Vir}{Vir} 
\DeclareMathOperator{\rep}{Rep} 
\DeclareMathOperator{\FP}{FP} 
\newcommand{\D}{\mathcal{D}} 
\newcommand{\I}{\mathfrak{i}} 
\newcommand{\cL}{\mathcal{L}_k}
\newcommand{\eL}{\mathcal{L}_k^{\even}}
\newcommand{\oL}{\mathcal{L}_k^{\odd}}
\newcommand{\cF}{\mathcal{F}_k}
\newcommand{\eF}{\mathcal{F}_k^{\even}}
\newcommand{\sF}{\mathcal{F}_k^{\odd}}
\newcommand{\cC}{\mathcal{C}_k}

\newcommand{\cD}{\mathcal{D}}
\newcommand{\voa}{vertex operator algebra}
\newcommand{\svoa}{vertex operator superalgebra}

\DeclareMathOperator{\id}{Id}

\title[$L_k\left(\osp(1 | 2)\right)$ from vertex tensor categories and Jacobi forms]{Representation theory of $L_k\left(\osp(1 | 2)\right)$ from vertex tensor categories and Jacobi forms }

\author{Thomas Creutzig, Jesse Frohlich and Shashank Kanade}
\address{Department of Mathematical and Statistical Sciences, University of Alberta,
	Edmonton, Alberta  T6G 2G1, Canada} 
\email{creutzig@ualberta.ca, jesse.frohlich@ualberta.ca, kanade@ualberta.ca}
\date{}

\allowdisplaybreaks

\begin{document}

\bibliographystyle{amsalpha}

\theoremstyle{plain}
\newtheorem*{introthm}{Theorem}
\newtheorem{obs}{Observation}
\newtheorem{thm}{Theorem}[section]
\newtheorem{prop}[thm]{Proposition}
\newtheorem{lem}[thm]{Lemma}
\newtheorem{cor}[thm]{Corollary}
\newtheorem{conj}[thm]{Conjecture}

\theoremstyle{definition}
\newtheorem{defi}[thm]{Definition}
\newtheorem{rem}[thm]{Remark}

\newcommand {\AL}{AL}
\newcommand {\ZZ}{\mathbb{Z}}
\newcommand {\RR}{\mathbb{R}}
\newcommand {\tr}{\text{tr}}
\newcommand {\sch}{\text{sch}}
\newcommand {\sltwo}{\mathfrak{sl}_2}
\newcommand {\cW}{\mathcal{W}}
\newcommand {\cX}{\mathcal{X}}
\newcommand {\cB}{\mathcal{B}}
\newcommand {\cS}{\mathcal{S}}
\newcommand {\cO}{\mathcal{O}_p}
\newcommand {\voas}{vertex operator algebras}
\newcommand {\svoas}{vertex operator superalgebras}
\newcommand{\Sing}{M(p)}
\newcommand{\Trip}{W(p)}
\newcommand{\voasl}{L_k(\mathfrak{sl}_2)}

\newcommand{\hopflink}{{\text{\textmarried}}}

\renewcommand{\baselinestretch}{1.2}

\date{}

\maketitle

\begin{abstract}
The purpose of this work is to illustrate in a family of interesting examples how to study the representation theory of  vertex operator superalgebras by combining 
the theory of vertex algebra extensions and modular forms. 

Let $L_k\left(\mathfrak{osp}(1 | 2)\right)$ be the simple affine vertex operator superalgebra of $\osp(1|2)$ at an admissible level $k$. We use a Jacobi form decomposition to see that this is a vertex operator superalgebra extension of $L_k(\mathfrak{sl}_2)\otimes \text{Vir}(p, (p+p')/2)$ where $k+3/2=p/(2p')$
and $\text{Vir}(u, v)$ denotes the regular Virasoro vertex operator algebra of central charge $c=1-6(u-v)^2/(uv)$. Especially, for a positive integer $k$, we get a regular vertex operator superalgebra and this case is studied further.

The interplay of the theory of vertex algebra extensions and modular data of the vertex operator subalgebra allows us to classify all simple local (untwisted) and Ramond twisted $L_k\left(\mathfrak{osp}(1 | 2)\right)$-modules and to obtain their super fusion rules. The latter are obtained in a second way from Verlinde's formula for vertex operator superalgebras. Finally, using again the theory of vertex algebra extensions, we find all simple modules and their fusion rules of the parafermionic coset $C_k = \text{Com}\left(V_L, L_k\left(\osp(1 | 2)\right)\right)$ where $V_L$ is the lattice vertex operator algebra of the lattice $L=\sqrt{2k}\Z$. 
\end{abstract}

\section{Introduction}

Understanding the representation theory of a given \voa{},
especially the tensor structure, is in general a difficult problem. 
However, if the given \voa{} is related to a known one by a standard construction then 
it is conceivable that much of the representation theory of the given \voa{} can be obtained from the known one. 
In a series of recent works \cite{CKL, CKLR, CKM} we have used vertex tensor categories to derive structural results about \voas{} and their extensions. The starting point is actually the important work of Huang, Kirillov and Lepowsky \cite{HKL} saying that \voa{} extensions of sufficiently nice \voas{} $V$ are in one-to-one correspondence with certain algebra objects $A$ in the representation category $\mathcal C$ of $V$. Moreover, the (untwisted or local) modules of the extended \voa{} are precisely objects in the category 
$\rep^0A$
of the local $A$-modules. This generalizes straightforwardly to extensions that are \svoas{} \cite{CKL}.
There is then an induction functor $\mathcal F$ (a tensor functor) that maps any object $X$ of the base category $\mathcal C$ to a not necessary local (super)algebra object $\mathcal F(X) \cong_{\mathcal C} A \boxtimes_{\mathcal C}X$, see \cite{KO}. In a recent work \cite{CKM}, the induction functor was studied from the \voa{} perspective. Most notably it was shown that $\rep^0A$ 
is braided equivalent to the category of extended \svoa{} modules and in addition the 
tensor product on $\rep^0A$ is exactly the $P(1)$-tensor product as defined in \cite{HLZ}
if mapping to local objects. 
The purpose of this work is to illustrate the usefulness of these recent results in very efficiently understanding the representation theories of two families of \svoas{} associated to $\osp(1|2)$.

\subsection{The affine \svoa{} $L_k\left(\osp(1 | 2)\right)$}

It is well-known that the simple affine \voa{} of a simple Lie algebra $\mathfrak{g}$ of level $k$ is regular if and only if $k$ is a positive integer \cite{FZ}. It is a presently unproven belief that the same statement holds for affine \svoa s  $L_k(\osp(1|2n))$ and that these are the only affine \svoa s that are regular. We note that fusion rules associated to $\widehat{\osp(1|2)}$ using coinvariants have been obtained in \cite{IK2}. This approach is believed but not known to give the \voa{} fusion rules.
 
Let $\cL$ be the simple affine \svoa{} of $\mathfrak{osp}(1|2)$ at level $k$ and $\eL$ its even component.
Characters of irreducible highest-weight modules of $\widehat{\osp(1|2)}$ for admissible level $k$ are known \cite{KW}. They converge in certain domains and allow for a meromorphic continuation to components of vector-valued Jacobi forms. It is then a Jacobi form decomposition problem to express these characters in terms of characters of $\voasl\otimes \Vir(p, \Delta/2)$. 
Here, $\Vir(u, v)$ denotes the simple and rational Virasoro \voa{} of central charge $1-6(u-v)^2/uv$ and the parameters are related to the level via $k+\frac{3}{2}=\frac{p}{2p'}$ and $\Delta=p+p'$. 
Our starting point is a character decomposition of simple highest-weight modules for $\widehat{\osp(1|2)}$, see Lemma \ref{lem:dec}. 
Since characters of non-isomorphic simple modules for both $\voasl$ and $M(p, \Delta/2)$ are linearly independent, we immediately get as a corollary that
\[
\text{Com}\left(\voasl,   \cL\right)  = \Vir(p, \Delta/2)\qquad \text{and}\qquad  \text{Com}\left(\Vir(p, \Delta/2),   \cL\right)  = \voasl.
\]
This result appeared in the physics literature \cite{CRS} and also in \cite{CL1}. 

Now let $k\in\ZZ_{>0}$.
It is easy to see that $\cL$ is $C_2$-cofinite, moreover, since an extension of a rational and $C_2$-cofinite \voa{} is rational, by \cite{DH, HKL, KO} another corollary is
that for $k$ in $\mathbb Z_{>0}$ the \svoa{} $\cL$ is rational and $C_2$-cofinite.
We then use the induction functors from $\voasl\otimes \Vir(p, \Delta/2)$ to $\eL$ and to $\cL$ to find a list of simple $\eL$ modules. 
There is then a powerful result \cite{DMNO} that relates the Frobenius-Perron dimension of the base category to the Frobenius-Perron dimension of the category of local modules for the algebra object. This dimension can be computed using the categorical $S$-matrix, that is the Hopf link, but this is related to the modular $S$-matrix due to Huang's theorem \cite{H1, H2}. We then use \cite{DMNO} together with the modular data of $\voasl\otimes \Vir(p, \Delta/2)$ to prove that we have found all simple $\eL$-modules. This immediately gives a complete list of simple local and Ramond twisted $\cL$-modules. All these results are obtained in Subsection \ref{sec:FP}.
The properties of the induction functor derived in \cite{CKM} then give us the fusion rules of simple $\cL$-modules (Subsection \ref{sec:fusion}) and modular $S$-matrix (Subsection \ref{sec:modular}). Alternatively, the fusion rules can be computed using Verlinde's formula for \svoas{} as derived in \cite{CKM} as a consequence of the Verlinde formula for regular vertex operator algebras due to Huang \cite{H1, H2}.

Extensions of these methods to admissible but non-integrable levels is
under investigation. For an interesting recent analysis of level $-5/4$
using methods different than ours see \cite{RSW}.

\subsection{The $L_k\left(\osp(1 | 2)\right)$ parafermion \voa}

The $L_k\left(\osp(1 | 2)\right)$ parafermion \voa{}s are called graded parafermions in physics \cite{CRS, FMW}.
Let $V$ be a regular vertex operator (super)algebra containing a lattice \voa{} $V_L$ of a positive definite lattice $L$ as a sub-vertex algebra. The coset Com$(V_L, V)$ is called the parafermion algebra of $V$. The best known case is if $V$ is a rational affine \voa. In that case $C_2$-cofiniteness \cite{ALY, DLY, DW} and rationality \cite{DR} have been understood a few years ago. If $V$ is the rational Bershadsky-Polyakov algebra, then these two properties have been proven indirectly \cite{ACL}. The very recent Theorem 4.12 of \cite{CKLR} settles this statement in generality and so we have that $C_k:=\text{Com}\left(V_L, L_k\left(\osp(1 | 2)\right)\right)$ is regular. Here, the lattice is $L=\sqrt{2k}\Z$. 
Section 4 of \cite{CKM} then allows us to deduce a complete list of inequivalent modules $C_{\lambda, r}$. These modules are parameterized by $\lambda$ in $L'/L$ and $r$ labelling the local $\cL$-modules. Their fusion rules are 
\[
C_{\nu, r} \boxtimes C_{\lambda, r'} \cong \bigoplus_{r'' =1 }^{2k+2} {N_{r, r'}}^{r''} C_{\lambda+\nu, r''},
\]
with $ {N_{r, r'}}^{r''}$ the fusion structure constants of $\cL$ derived earlier. One can also express characters of the $C_{\lambda, r}$ and their modular transformations immediately in terms of those of $\cL$. 
All these results are specific instances of results in Section 4 of \cite{CKM} and are presented in Section \ref{sec:para}.

\subsection*{Acknowledgments}
T.C.\ is supported by the Natural Sciences and Engineering Research Council of Canada (RES0020460).
S.K.\ is supported in part by PIMS post-doctoral fellowship and by Endeavour Research Fellowship (2017)
awarded by the Department of Education and Training, Australian Government. 
We thank Robert McRae and David Ridout for illuminating discussions.

\section{Decomposition of Jacobi forms associated to $\osp(1 | 2)$}

Let us first recall the necessary known data of the relevant \voa s and Lie super algebras from the literature. 

\subsection{Characters of affine Lie super algebra $\widehat{\osp(1|2)}$ modules}
\label{subsec:charosp12}

The conformal field theory pendant of the affine \svoa{} of $\widehat{\osp(1|2)}$ has appeared in \cite{ENO, ER}.
It will turn out that all we need in order to understand much of the representation theory of admissible integer level $L_k\left(\osp(1|2)\right)$ is a certain decomposition of the vacuum module character. Strong reconstruction theorem gives that the simple quotient of the level $k$ vacuum Verma module of an affine Lie super algebra carries the structure of a simple \svoa{} \cite{FBZ}.
Admissible levels are those where the vacuum Verma module itself is not simple and the character of its simple quotient can be meromorphically continued to a meromorphic Jacobi form. In the case of $\widehat{\osp(1|2)}$ these are
\begin{equation}\label{eqn:kp}
  k+\frac 32=\frac{p}{2p'}
\end{equation}
where $p>1$ and $p'$ are positive coprime integers with $p+p'\in 2\Z$ and $p$,
$\frac{p+p'}{2}$ coprime. The highest-weight representations $M_{j_{r,s}}$ are indexed by isospins $j_{r,s}$,
with
\begin{equation*}
    1\le r \le p-1, \qquad
    0\le s \le p'-1, \qquad \text{and} \qquad
    r+s \in 2\Z + 1.
\end{equation*}
From \cite{ER, KW}, the characters of these representations are given by
\begin{equation} \label{eq:char}
    \ch[M_{j_{r,s}}]=\frac{
        \Theta_{b_+,a}\left(\dfrac{z}{2p'},\dfrac{\tau}{2}\right)
        -\Theta_{b_-,a}\left(\dfrac{z}{2p'},\dfrac{\tau}{2}\right)
    }{\Pi(z,t)}
\end{equation}
with standard Jacobi theta functions
\[
    \Theta_{r,s}(z,\tau)
    =\sum_{m\in\Z}
        w^{s(m+\frac{r}{2s})}
        q^{s(m+\frac{r}{2s})^2}
        \]
        and the Weyl super denominator
\[
        \Pi(z,\tau)=\Theta_{1,3}\left(\frac z2,\frac \tau 2 \right)
            -\Theta_{-1,3}\left(\frac z2,\frac \tau 2 \right)
            = w^{\frac{1}{4}} q^{\frac{1}{24}} \prod_{n=1}^{\infty} \frac{
            \left(1-q^{n}\right) \left(1-wq^{n}\right)
            \left(1-{w}^{-1}q^{n-1}\right) }{
            \left(1+w^{\frac 12}q^n\right)
            \left(1+w^{-\frac 12}q^{n-1}\right)}
\]
\[
\text{with} \qquad \qquad\qquad
    w=e^{2\pi i z}, \qquad q=e^{2\pi i \tau}, \qquad b_\pm=\pm p'r-ps \qquad\text{and} \qquad a=pp'. \qquad\qquad\qquad\qquad\qquad{}
\]
These characters are analytic functions in the domain $1 \le |w| \le |q^{-1}|$ and can be meromorphically continued to meromorphic Jacobi forms.  Their modular transformations are a straight-forward computation, we will obtain them in the positive integer level case directly from vertex tensor category theory. 

\subsection{The rational Virasoro \voa{} $\Vir(p, u)$}

A good reference here is \cite{IK}.
Let $u, p\in \mathbb Z_{\geq 2}$ coprime. Then, the simple Virasoro \voa{} at central charge 
\[
c=1 - 6 \frac{(u-p)^2}{up}
\]
is regular \cite{Wang}. Simple modules are denoted by $V_{r, s}$ for $1\leq r \leq u-1$ and $1\leq s\leq p-1$ and one has the relations $V_{r, s} \cong V_{u-r, p-s}$. We denote the set of inequivalent module labels by $I_{u, p}$.
We define numbers by 
\[
{N^{w \ \, t''}_{t, t'}} = \begin{cases} 1 & \quad \text{if} \ |t-t'| +1 \leq t'' \leq \text{min} \{ t+t'-1, 2w-t-t'\}, t+t'+t'' \ \text{odd} \\
0 & \quad \text{otherwise}
 \end{cases}
\]
that allow to express the fusion rules as follows
\[
V_{r, s} \boxtimes_{\Vir} V_{r', s'} \cong  \bigoplus_{r''=1}^{u-1} \bigoplus_{s''=1}^{p-1} {N^{u \ \, r''}_{r, r'}}{N^{p \ \, s''}_{s, s'}} V_{r'', s''}.
 \]
With $p,p'$ as in Subsection \ref{subsec:charosp12} set 
 $\Delta=p+p'$ and $u=\Delta/2$. 
Then characters are \cite{IK}
\begin{equation}
\begin{split}
    \ch[V_{r, s}] &=  \frac{\left(
        \Theta_{2pr-\Delta s,\Delta p}-
        \Theta_{-2pr-\Delta s,\Delta p}
    \right)\left(0,\frac{\tau}{2}\right)}{\eta(q)}
\end{split}
\end{equation}
and their modular $S$-transformation is
\begin{equation}
\ch[V_{r, s}]\left(-\frac{1}{\tau}\right) = \sum_{(r', s') \in I_{u, p}} S^\chi_{(r, s), (r', s')} \ch[V_{r', s'}](\tau)
\end{equation}
 with modular $S$-matrix entries
\begin{equation}
\label{eqn:virS}
 S^\chi_{(r, s), (r', s')} = -2 \sqrt{\frac{2}{up}} (-1)^{rs'+sr'} \sin\left(\frac{\pi p}{u} rr'  \right) \sin\left(\frac{\pi u}{p} ss'  \right). 
\end{equation}
The conformal weights of $V_{r,s}$ are given by
 $ h_{r,s} = \frac{(us-pr)^2 - (u-p)^2}{4up}$.

 \subsection{Characters of admissible level $L_k(\mathfrak{sl}_2)$ modules}
 
 We use \cite{KW, CR2} as reference. Let the level $k$ be as in \eqref{eqn:kp}, i.e., $k$ satisfies  
 \[
 k+2= k+\frac{3}{2}+\frac{1}{2}= \frac{p}{2p'}+\frac{1}{2} = \frac{p+p'}{2p'}=\frac{u}{p'}.
 \]
 The highest-weight representations of $L_k(\mathfrak{sl}_2)$ are $\cD^+_{r, s}$ with $1\leq r\leq u-1$ and $0\leq s\leq p'-1$
where $\cD^+_{r, 0};=\mathcal L_{r,0}$ deserve a special name as they have finite-dimensional grade zero subspace.
 Let  
\begin{align*}
    \vartheta_2(z,\tau)&=\sum_{n\in\Z}
    w^{\left(n+\frac 12\right)}
    q^{\frac 12\left(n+\frac 12\right)^2}, \qquad \vartheta_1(z, \tau) = -\vartheta_2\left(z+\frac{1}{2}, \tau\right)
\end{align*}
be the standard theta functions, 
then characters are
\[
\ch[\cD^+_{r, s}](z, \tau) =   \frac{\left(
        \Theta_{2p'r-\Delta s,\Delta p'}-
        \Theta_{-2p'r-\Delta s,\Delta p'}
    \right)\left(\dfrac{z}{2p'},\dfrac{\tau}{2}\right)}{\I \vartheta_1(w,q)}.
\]
If $k\in\mathbb{Z}_{>0}$, then the $\cD^+_{r, 0}:=\mathcal L_{r,0}$ with $1\leq r\leq k+1$ of conformal dimensions $\dfrac{r^2-1}{4(k+2)}$ exhaust inequivalent simple modules up to isomorphisms and their modular $S$-transformations are given by
\begin{equation}
\ch[\mathcal L_{r,0}]\left(0, -\frac{1}{\tau}\right) = \sum_{r'=1}^{k+1} S^{\chi,\mathfrak{sl}_2}_{r, r'} \ch[\mathcal L_{r',0}](0, \tau)
\qquad\text{with}
\qquad
S^{\chi,\mathfrak{sl}_2}_{r, r'}  = \sqrt{\frac{2}{k+2}} \sin\left(\frac{\pi}{k+2}rr'\right). 
\label{eqn:sl2S}
\end{equation}
\begin{lem}\label{lem:dec}
For $k$ as in \eqref{eqn:kp}, we have the following character decomposition
\begin{align*}
    \ch[M_{j_{r,s}}](w, q)=\sum_{i=1}^{\frac\Delta2-1}
    \ch[\D^+_{i,s}](w^{\frac 12},q)
    \ch[\Vir^{\left(p,\frac\Delta2\right)}_{i,r}](q)
\end{align*}
\end{lem}
\begin{proof}
We compare left and right-hand side each multiplied by $\vartheta_2\left(\frac z2,\tau\right)$.
\begin{align}
   \Theta_{b_\pm,a}\left(\frac{z}{2p'},\frac\tau2\right)\vartheta_2\left(\frac z2,\tau\right)
   &=\sum_{m,n\in\Z}
        w^{\frac{a}{2p'}\left(m+\frac{b_\pm}{2a}\right)+\frac 12\left(n+\frac12\right)}
        q^{\frac{a}{2}\left(m+\frac{b_\pm}{2a}\right)^2+\frac 12\left(n+\frac12\right)^2} \nonumber\\ 
    &=\sum_{m,n\in\Z}
        w^{\frac{p}{2}\left(m+\frac{b_\pm}{2a}\right)+\frac 12\left(n+\frac 12\right)}
        q^{\frac{pp'}{2}\left(m+\frac{b_\pm}{2pp'}\right)^2+\frac 12\left(n+\frac12\right)^2}\nonumber\\ 
    &=\sum_{m,n\in\Z} w^{\frac{1}{2}\left(pm+\frac{b_\pm+p'}{2p'}+n\right)}
        q^{\frac{1}{2}\left(pp'm^2+m b_\pm+\frac{b_\pm^2}{4pp'}+n^2+n+\frac14\right)}.
         \label{eq:osp}
\end{align}
At this point, we make a change of variables. Let $x=pm+n$. We wish to find
$y=cm+dn$ such that $pp'm^2+n^2$ has no mixed terms in $x$ and $y$. Let
$\Delta=c-dp$ be the determinant of this change of variables. Then
$m=\Delta^{-1}(y-dx)$ and $n=\Delta^{-1}(cx-py)$. From this, we see:
\begin{align} \label{eq:exp}
    pp'm^2+n^2
    &=\frac{\left(pp'd^2+ c^2\right)x^2+p\left(p'+p\right)y^2-p\left(c+p'd\right)xy}{\Delta^{2}}
\end{align}
So, we see that $c+p'd$ needs to vanish. Choosing $c=p'$ and $d=-1$, this condition is
fulfilled. Therefore, $\Delta=p+p'$. The change of variables becomes
$m=\Delta^{-1}(x+y)$ and $n=\Delta^{-1}(p'x-py)$. With the observations that
$p'=-p\pmod\Delta$ and $p$ is relatively prime to $\Delta$, we see that this change
of variables is an invertible map $\varphi:\Z^2\to L=\left\{(x,y)\in\ZZ^2\middle|x+y\in\Delta\Z\right\}$.
So, one may simply sum over $L$ in $x$ and $y$. Simplifying \eqref{eq:exp}, we
now have:
\begin{align*}
    pp'm^2+n^2
    &=\frac{p'}{\Delta}x^2+\frac{p}{\Delta}y^2
\end{align*}
From this, one easily obtains the exponent for $q$:
\begin{align*}
    pp'm^2+n^2+mb_\pm+\frac{b_\pm^2+pp'}{4pp'}+n
        &=\frac{p'}{\Delta}\left(x+\frac{b_\pm+p'}{2p'}\right)^2
        +\frac{p}{\Delta}\left(y+\frac{b_\pm-p}{2p}\right)^2.
\end{align*}
Returning to \eqref{eq:osp}, and letting $x=\Delta s+i$ and $y=\Delta t-i$, for
$i=1,\dots,\Delta$,
\begin{align} \label{eq:final}
\Theta_{b_\pm,a}&\left(\frac{z}{2p'},\frac{\tau}{2}\right)\vartheta_2\left(\frac z2,\tau\right)
       =\sum_{(x,y)\in L}
        w^{\frac{1}{2}\left(x+\frac{b_\pm+p'}{2p'}\right)}
        q^{\frac{p'}{2\Delta}\left(x+\frac{b_\pm+p'}{2p'}\right)^2}
        q^{ \frac{p}{2\Delta}\left(y+ \frac{b_\pm-p}{2p}\right)^2} \nonumber \\
    &=\sum_{i=1}^{\Delta}\left(\sum_{s\in\Z}
        w^{\frac{\Delta p'}{2p'}\left(s+\frac{b_\pm+p'(2i+1)}{2\Delta p'}\right)}
        q^{\frac{\Delta p'}{2}\left(s+\frac{b_\pm+p'(2i+1)}{2\Delta p'}\right)^2}\right)
        \left(\sum_{t\in\Z}
        q^{\frac{\Delta p}{2}\left(t+\frac{b_\pm-p(2i+1)}{2\Delta p}\right)^2}\right)
        \nonumber \\
    &=\sum_{i=1}^{\Delta}
        \Theta_{b_\pm+p'(2i+1),\Delta p'}\left(\frac{z}{2p'},\frac{\tau}{2}\right)
        \Theta_{b_\pm-p(2i+1),\Delta p}\left(0,\frac{\tau}{2}\right) \nonumber \\
    &=\sum_{i=1}^{\Delta}
        \Theta_{-ps+p'(2i\pm r+1),\Delta p'}\left(\frac{z}{2p'},\frac{\tau}{2}\right)
        \Theta_{\pm p'r-p(2i+s+1),\Delta p}\left(0,\frac{\tau}{2}\right)
\end{align}
Note that the theta-functions in \eqref{eq:final} are $\Delta$-periodic in $i$,
so the whole sum in \eqref{eq:final} is $1$-periodic in $i$.
Next, we want to consider the difference between the two results, according to
\eqref{eq:char}.
\allowdisplaybreaks
\begin{equation}\nonumber
\begin{split}
 \text{Let}& \   X:=  \left(\Theta_{b_+,a}\left(\frac{z}{2p'},\frac{\tau}{2}\right)-
    \Theta_{b_-,a}\left(\frac{z}{2p'},\frac{\tau}{2}\right)\right)
    \vartheta_2\left(\frac{z}{2},\tau\right) \ \text{which expands to} \qquad\qquad\qquad \\
    X&=\sum_{i=1}^{\Delta}
        \Theta_{-ps+p'(2i+r+1),\Delta p'}\left(\frac{z}{2p'},\frac{\tau}{2}\right)
        \Theta_{p'r-p(2i+s+1),\Delta p}\left(0,\frac{\tau}{2}\right) \\
& \qquad\qquad\qquad\qquad    -\sum_{i=1}^{\Delta}
        \Theta_{-ps+p'(2i-r+1),\Delta p'}\left(\frac{z}{2p'},\frac{\tau}{2}\right)
        \Theta_{-p'r-p(2i+s+1),\Delta p}\left(0,\frac{\tau}{2}\right).
\end{split}
\end{equation}
After applying some symmetries to the second term, such as negation of the
second factor's first index and mapping $i \mapsto -i-s-1$:
\begin{align}\nonumber
  X&=  \sum_{i=1}^{\Delta}
    \left(
        \Theta_{-ps+p'(2i+r+1),\Delta p'}\left(\frac{z}{2p'},\frac{\tau}{2}\right)-
        \Theta_{-ps-p'(2i+r+2s+1),\Delta p'}\left(\frac{z}{2p'},\frac{\tau}{2}\right)
    \right)\nonumber\\
    &\quad\quad\cdot \Theta_{p'r-p(2i+s+1),\Delta p}\left(0,\frac{\tau}{2}\right)\\
  &=  \sum_{i=1}^{\Delta}	
    \left(
        \Theta_{-\Delta s+p'(2i+r+s+1),\Delta p'}\left(\frac{z}{2p'},\frac{\tau}{2}\right)-
        \Theta_{-\Delta s-p'(2i+r+s+1),\Delta p'}\left(\frac{z}{2p'},\frac{\tau}{2}\right)
    \right)\nonumber\\
    &\quad\quad\cdot\Theta_{p'r-p(2i+s+1),\Delta p}\left(0,\frac{\tau}{2}\right).
\end{align}
By assumption, $r+s\in2\Z+1$. Thus, letting $t = \frac{r+s+1}{2}$, we have
\begin{align}\nonumber
    X&=\sum_{i=1}^{\Delta}
    \left(
        \Theta_{-\Delta s+2p'(i+t),\Delta p'}\left(\frac{z}{2p'},\frac{\tau}{2}\right)-
        \Theta_{-\Delta s-2p'(i+t),\Delta p'}\left(\frac{z}{2p'},\frac{\tau}{2}\right)
    \right)\Theta_{p'r-p(2i+s+1),\Delta p}\left(0,\frac{\tau}{2}\right)\\
&=    \sum_{i=1}^{\Delta}
    \left(
        \Theta_{2p'i-\Delta s,\Delta p'}\left(\frac{z}{2p'},\frac{\tau}{2}\right)-
        \Theta_{-2p'i-\Delta s,\Delta p'}\left(\frac{z}{2p'},\frac{\tau}{2}\right)
    \right)\Theta_{\Delta r-2pi,\Delta p}\left(0,\frac{\tau}{2}\right).
\end{align}
For $i=\Delta, \frac\Delta2$, the first factor
vanishes. We can combine summands to obtain
\begin{align*}
   X&= \sum_{i=1}^{\frac\Delta2-1}
    \left(
        \Theta_{2p'i-\Delta s,\Delta p'}\left(\frac{z}{2p'},\frac{\tau}{2}\right)-
        \Theta_{-2p'i-\Delta s,\Delta p'}\left(\frac{z}{2p'},\frac{\tau}{2}\right)
    \right) \Theta_{\Delta r-2pi,\Delta p}\left(0,\frac{\tau}{2}\right) \\
   &\quad -\sum_{i=1}^{\frac\Delta2-1}
    \left(
        \Theta_{2p'i-\Delta s,\Delta p'}\left(\frac{z}{2p'},\frac{\tau}{2}\right)-
        \Theta_{-2p'i-\Delta s,\Delta p'}\left(\frac{z}{2p'},\frac{\tau}{2}\right)
    \right) \Theta_{\Delta r+2pi,\Delta p}\left(0,\frac{\tau}{2}\right)\\
&=\sum_{i=1}^{\frac\Delta2-1}
    \left(
        \Theta_{2p'i-\Delta s,\Delta p'}-
        \Theta_{-2p'i-\Delta s,\Delta p'}
    \right) \left(\frac{z}{2p'},\frac{\tau}{2}\right)
    \left(
        \Theta_{2pi-\Delta r,\Delta p}-
        \Theta_{-2pi-\Delta r,\Delta p}
    \right)\left(0,\frac{\tau}{2}\right).
\end{align*}
Comparing with the $L_k(\mathfrak{sl}_2)$ and Virasoro characters now easily gives the claim.
\end{proof}
Since characters of modules in the category $\mathcal O$ of $L_k(\mathfrak{sl}_2)$ at admissible level and those of the rational Virasoro vertex algebra both determine representations, we
have the analogous result on the level of modules
\begin{equation} \label{eq:decomp}
    M_{j_{r,s}}=\bigoplus_{i=1}^{\frac\Delta2-1}
    \D^+_{i,s} \otimes V_{i,r}.
\end{equation}
Here we remark that in general, characters of simple \voa{} modules are not necessarily linearly independent and in the case 
$L_k(\mathfrak{sl}_2)$ for $k\in\QQ\backslash\ZZ$
this question is subtle. In that case there are simple modules whose characters only converge to an analytic function in a certain domain (depending on the module). These characters can then be meromorphically continued to meromorphic Jacobi forms and many different module characters will have same meromorphic continuation. These modules are then distinguished by the domain in which the character converges to an analytic function. For each meromorphic Jacobi form in question there is a unique simple module whose character converges in the domain $1\le |w|\le |q^{-1}|$ and coincides in that domain with the meromorphic Jacobi form. These interesting subtleties are discussed thoroughly in \cite{CR1, CR2}.

\section{Representation theory of $\cL$ for $k\in\mathbb Z_{>0}$}

Let $k\in\mathbb Z_{>0}$ and consider 
\begin{equation}
   \cL=\bigoplus_{i=1}^{k+1}
    \mathcal L_{i,0} \otimes V_{i,1}
\qquad
\text{and} 
\qquad
   \eL=\bigoplus_{\substack{ i=1 \\ i \ \text{odd} }}^{k+1}
    \mathcal L_{i,0} \otimes V_{i,1}.
\end{equation}
$\cL$ can be endowed with the structure of 
$L_k(\osp(1|2))$
and $\eL$ as its even 
subalgebra
is simple as well.  
It follows from \cite{CarM, CKLR} that
\begin{equation}
   \oL=\bigoplus_{\substack{ i=1 \\ i \ {\even} }}^{k+1}
    \mathcal L_{i,0} \otimes V_{i,1}
\end{equation}
is an order two simple current of $\eL$. 
Let $\cC$ be the modular tensor category of  the \voa{} $\voasl\otimes \Vir(p, \Delta/2)$.
We can view $\cL$ as a superalgebra object in both $\cC$ and in the vertex tensor category of $\eL$ while $\eL$ is also an algebra object in $\cC$.
Denote the corresponding induction functors by $\cF, \sF, \eF$.  Our aim is to use the theory of vertex algebra extensions \cite{CKM} to understand the representations of $\cL$. 

Recall that the categorical twist $\theta$ acts on \voa{} modules as $e^{2\pi i L(0)}$ with $L(0)$ the Virasoro zero-mode that provides the grading by conformal weight of modules. 
Conformal dimension of $\mathcal{L}_{i,0}\otimes V_{i,r}$ is
\begin{equation}
\label{eqn:confdimLioVir}\dfrac{1}{4}\left(2i^2-2ir +\dfrac{(k+2)(r^2-1)}{2k+3}\right) \qquad  (\text{mod}\  \ZZ).
\end{equation}
We have that $\theta_{\eL}=\id_{\eL}$.
Also, the calculation \eqref{eqn:dimLkeven} below shows that ``dimension'' $\dim_{\cC}(\eL)$ is non-zero.
Using \cite[Lem.\ 1.20, Thm.\ 4.5]{KO} and \cite[Thm.\ 3.65]{CKM} 
the representation category of untwisted modules of $\eL$ is modular.

Now, with $1\leq r\leq 2k+2$, we have the following (possibly non-local) modules for $\cL$ via induction
\[
M_r = \cF\left( \mathcal L_{1,0} \otimes V_{1,r} \right) = \bigoplus_{i=1}^{k+1}
    \mathcal L_{i,0} \otimes V_{i,r} = M_r^{{\even}} \oplus M_r^{{\odd}}
\]
with the $\eL$-modules 
\begin{align*}
M_r^{{\even}}  &= \eF\left( \mathcal L_{1,0} \otimes V_{1,r} \right) = \bigoplus_{\substack{ i=1\\ i \ {\odd}} }^{k+1}
    \mathcal L_{i,0} \otimes V_{i,r},\quad\quad
    M_r^{{\odd}}= 
     \bigoplus_{\substack{ i=1\\ i \ {\even}} }^{k+1}
    \mathcal L_{i,0} \otimes V_{i,r}.
\end{align*}
All $M_r$ and $M_r^{{\even}}$ 
are simple $\cL$-, respectively, $\eL$-modules by \cite[Prop.\ 4.4]{CKM}. 
\cite[Lem.\ 4.26]{CKM} now implies that $M_r^{{\odd}}\cong \oL\boxtimes_{\eL}M_r^{{\even}}$,
in particular, they are simple. We now have that $\cF\left( M_r^{{\even}}  \right) \cong M_r$.

$M_r$ is local if it is graded by a coset of integers and twisted 
if both $M_r^{{\even}}$ and $M_r^{{\odd}}$ are each graded by a single but different cosets of integers. 
From \eqref{eqn:confdimLioVir} we have that $M_r$ is local if and only if $r$ is odd and otherwise it is twisted. 
We aim to prove that $M_r$ exhaust all local/twisted modules (up to isomorphisms) of $\cL$ and to determine their modular data and fusion rules. 

Using the fact that induction is a tensor functor \cite[Thm.\ 3.65, Thm.\ 3.68]{CKM} we have that:
\begin{align}
&M_r^{\even}\boxtimes_{{\eL}}M_{r'}^{\even}=\eF(\mathcal L_{1,0}\otimes V_{1,r})\boxtimes_{{\eL}}\eF(\mathcal L_{1,0}\otimes V_{1,r'})\\
&=\eF((\mathcal L_{1,0}\otimes V_{1,r})\boxtimes (\mathcal L_{1,0}\otimes V_{1,r'}))\nonumber
=\eF(\bigoplus_{r''=1}^{k}N_{r,r'}^{k+1\,\,r''}\mathcal L_{1,0}\otimes V_{1,r''})
=\bigoplus_{r''=1}^{k}N_{r,r'}^{r''}M_{r''}^{\even},
\end{align}
and similarly,
\begin{align}
M_r\boxtimes_{{\cL}}M_{r'}&=\bigoplus_{r''=1}^{k}N_{r,r'}^{r''}M_{r''}.
\label{eqn:Mrfus}
\end{align}
In particular, the abelian (super)categories generated by $M_r$ and $M_r^{\even}$ 
are closed under taking tensor products in respective categories.
Moreover, since $M_r^{\odd}\cong \oL\boxtimes_{\eL}M_r^{\even}$, 
the abelian category generated by $M_r^{\even/\odd}$ also has the same property
and moreover, it is closed under duals:
\begin{align*}
(M_r^{\even})^*&=\eF(\mathcal L_{1,0}\otimes V_{1,r})^*=\eF((\mathcal L_{1,0}\otimes V_{1,r})^*)
=\eF(\mathcal L_{1,0}\otimes V_{1,r})=M_r^{\even}\\
(M_r^{\odd})^*&=(\oL\boxtimes_{\eL} M_r^{\even})^*\cong {\oL}^*\boxtimes_{\eL}(M_r^{\even})^*\cong M_r^{\odd},
\end{align*}
where the second equality in the first equation follows by \cite[Prop.\ 2.17]{CKM}, and for the 
second equation we use the fact that $\eL$ is an order two simple current.
In effect, this category is a fusion category. We denote this category by $\cS_k$ for future use.

\subsection{Frobenius-Perron dimension of $\mathcal L^{\even}$}\label{sec:FP}

We now use modular properties of characters to classify all simple inequivalent modules of $\eL$ and of $\cL$ as well.

The categorical dimension of a \voa{} module $X$ is the Hopf link $S^{\hopflink}_{\mathbf{1},X}$ which coincides with the modular $S$-matrix expression
up to a constant as $S^{\hopflink}_{\mathbf{1},X}=\frac{S^{\chi}_{\mathbf{1},X}}{S^{\chi}_{\mathbf{1}, \mathbf{1}}}$. Here $\mathbf{1}$ stands for the tensor unit, i.e., the \voa{} itself. 
Using \eqref{eqn:virS} and \eqref{eqn:sl2S}, one has
\begin{equation}\label{eqn:dimLkeven}
\dim_{\cC} \eL  =\frac1{\sin^2\left(\frac{\pi}{k+2}\right)}
\sum_{\substack{l=1 \\ l\in2\Z+1}}^{k+1}
\sin^2\left(\frac{\pi l}{k+2}\right),
\end{equation}
which we can simplify using the following lemma.
\begin{lem}\label{lem:sin2}
	The following holds:
	\begin{equation}
	\sum_{\substack{l=1 \\ l\in2\Z+1}}^{k+1} \sin^2\left(\frac{\pi l}{k+2}\right)
	= \frac{k+2}{4}\qquad \text{and} \qquad \sum_{{l=1}}^{k+1} \sin^2\left(\frac{\pi l}{k+2}\right)
	= \frac{k+2}{2}.
	\end{equation}
\end{lem}
\begin{proof} 
	The second equation is 
	\begin{align*}
	\sum_{l=1}^{k+1}\sin^2\left(\frac{\pi l}{k+2}\right)
	=-\frac{1}{4}\sum_{l=1}^{k+1} \left(-2
	+e^{2\pi\I\frac{l}{k+2}}
	+e^{-2\pi\I\frac{l}{k+2}}\right) = \frac{k+1}{2} + \frac 14 + \frac 14 = \frac{k+2}{2}.
	\end{align*}
	If $k$ is even, then the first equation follows immediately from the second one. 
	For odd $k$ we compute
	\begin{equation} \nonumber
	\begin{split}
	\sum_{\substack{l=1 \\ l\in2\Z+1}}^{k+1} \sin^2\left(\frac{\pi l}{k+2}\right)
	&=\sum_{r=0}^{\frac{k-1}2} \sin^2\left(\frac{\pi (2r+1)}{k+2}\right)=\frac{k+1}{4}-
	\frac{1}{4}\sum_{r=0}^{\frac{k-1}2}
	\left(e^{\frac{2\pi\I(2r+1)}{k+2}}+ e^{-\frac{2\pi\I(2r+1)}{k+2}}\right)
	\end{split}
	\end{equation}
	and 
	\begin{equation} \nonumber
	\begin{split}
	\sum_{r=0}^{\frac{k-1}2}
	\left(e^{\frac{2\pi\I(2r+1)}{k+2}}+ e^{-\frac{2\pi\I(2r+1)}{k+2}}\right) &= 
	2 \sum_{r=0}^{k+1} e^{\frac{2\pi\I r}{k+2}}
	-\sum_{r=0}^{\frac{k+1}2} \left(e^{\frac{2\pi\I(2r)}{k+2}}+e^{-\frac{2\pi\I(2r)}{k+2}}\right) 
	= -1.\qquad\qquad  \qedhere
	\end{split}
	\end{equation}
\end{proof}

The Frobenius-Perron dimension (FP for short) of an object is uniquely characterized by the following.
For simple modules $X$, $\FP(X)\in\RR_{>0}$ such that the map $W\mapsto$ $\FP(W)$ is a one-dimensional representation of the tensor ring, see \cite{ENO}.
The Frobenius-Perron dimension of a modular tensor category $\mathcal{C}$ is
\[
\FP(\mathcal C) =\sum  \FP(X)^2
\]
where the sum is over all inequivalent simple objects of $\mathcal{C}$. 
For \voa s $V$, this dimension is easy to find provided there exists a unique simple module $Z$ such that it has strictly lowest conformal dimension among all simple modules for $V$ and provided the following limits exist \cite{DJX}. In that case for any simple object $X$ one has
\[
\text{adim}(X) := \lim_{\tau\rightarrow 0^+}\frac{\ch[X](\tau)}{\ch[V](\tau)} =  \lim_{\tau\rightarrow -\infty}\frac{\ch[X](-1/\tau)}{\ch[V](-1/\tau)} = \frac{S^\chi_{X, Z}}{S^\chi_{V, Z}}
\]
and $\text{adim}(X) =  \FP(X)$. 

Recall that the $M_r^{\text{even/odd}}$ form a sub-tensor category, called $\cS_k$, of the full category of local modules of $\eL$. 
We will compute the $\FP(\cS_k)$ and show that it agrees with $\FP(\rep^0\eL)$ so that the two categories must coincide. 
We also recall:
\begin{thm}\textup{\cite[Cor.\,3.30]{DMNO}}\label{lem:FPrepA}
Let $A$ be a connected (sometimes also called haploid) \'etale (that is, commutative and separable) 
algebra in a modular tensor category $\mathcal C$, 
with $\rep^0A$ denoting the category of its local (also known as untwisted or dyslectic) modules in $\mathcal C$,
then 
\[
\FP\left(\rep^0(A)\right) = \frac{\FP(\mathcal C)}{\FP_\mathcal C(A)^2}.
\]
\end{thm}

\begin{lem}\label{lem:FPCk}
The Frobenius-Perron dimensions of $\cC$ and $\eL$ are
\[
\FP(\cC) = \frac{(k+2)^2(2k+3)}{16\sin^4\left(\frac{\pi}{k+2}\right)\sin^2\left(\frac{\pi}{2k+3}\right) } \qquad\text{and} \qquad \FP_{\cC}\left( \eL \right)  =\frac{k+2}{4\sin^2\left(\frac{\pi}{k+2}\right)}.\]
\end{lem}
\begin{proof}
The simple module of lowest conformal weight for $\voasl$ is $\voasl$ itself as this is a unitary \voa.
For the Virasoro algebra, we will prove that $h_{1, 2}\leq h_{r, s}$ and equality holds if and only if $V_{r, s}\cong V_{1, 2}$.
The conformal weights of $V_{r,s}$ are given by
\begin{equation*}
  h_{r,s} = \frac{\left((k+2)s-(2k+3)r\right)^2 - \left(2k+3 - (k+2)\right)^2}{4(k+2)(2k+3)}
\end{equation*}
where $1 \le r
\le k+1 =: t$ and $1 \le s \le 2k+2 = 2t$. Minimizing $h_{r,s}$ amounts to minimizing $X_{r,s} := |(k+2)s-(2k+3)r| = |t(s-2r)+(s-r)|$.
\begin{enumerate}
\item[Case 1.] $2r-s=0$:  Since $X_{r,s} = r$, we see $X_{1,2} = 1$.
\item[Case 2.] $2r-s=1$:  $X_{r,2r-1} = |r-1-t|
    \leadsto X_{t,2t-1} = 1$. Note that $V_{r,s} \cong V_{t-r+1,2t-s+1}$, so $V_{1,2} \cong V_{t,2t-1}$.
\item[Case 3.] $2r-s=n\ge2$:  $X_{r,2r-n} = |n(t+1)-r| \ge (n-1)t+n \ge t+2 > 2$.
\item[Case 4.] $2r-s=-n<0$: $X_{r,2r+n} = |t(n+1)+r| > 1$.
\end{enumerate}
The claim now follows directly with the previous Lemma together with the modular $S$-matrices of $L_k(\mathfrak{sl}_2)$ and the Virasoro algebra modules. 
\end{proof}
We have already shown that $\cS_k$, the full abelian sub-category formed by $M_r^{\text{even/odd}}$-modules
inside $\eL$-modules is a fusion category.
Considering the homomorphism of Grothendieck rings 
induced by the inclusion of $\cS_k\hookrightarrow \rep^0\eL$,
from \cite[Prop. 3.3.13(i)]{EGNO},  we now immediately have that
$\FP$-dims of simples in $\mathcal{S}_k$ are equal to their $\FP$-dims as objects of $\rep^0\eL$
and hence that $\FP(\cS_k)\leq \FP(\rep^0\eL)$ with equality iff the two categories
are equal.

Now, with the $S$-matrix calculations from \eqref{eqn:SMeMe} and \eqref{eqn:SMMother},
it is clear that $X\mapsto \dfrac{S^{\chi}_{X,M_2^{{\even}}}}{S^{\chi}_{\eL,M_2^{{\even}}}}\in\RR_{>0}$ 
whenever $X=M_r^{\text{even/odd}}$, and 
moreover, this map preserves tensor products by properties of $S$ matrices
for $\eL$. Therefore, by uniqueness of $\FP$ dimensions, this 
map precisely gives us $\FP$ dimensions of modules in $\cS_k$
(which are equal to their $\FP$-dims considered as modules for $\eL$).

The following Lemma is verified in a very similar manner as the previous one, using the modular $S$-matrix derived in the next section. 
\begin{lem}\label{lem:FPMrevenodd}
We have the following:
\[
\FP(\cS_k)=\frac{2k+3}{\sin^2\left(\frac{\pi}{2k+3} \right)}.
\]
\end{lem}
\begin{proof} From the $S$-matrices in the next section \eqref{eqn:SLeven}, \eqref{eqn:SMeMe} and
	\eqref{eqn:SMMother} and Lemma \ref{lem:sin2}, we have that
\begin{align*}
\FP&(\cS_k) =\sum_{1\leq r \leq 2k+2} 
\left(\dfrac{S^{\chi}_{M_r^{{\even}},M_2^{{\even}}}}{S^{\chi}_{M_1^{{\even}},M_2^{{\even}}}}\right)^2
+\left(\dfrac{S^{\chi}_{M_r^{{\odd}},M_2^{{\even}}}}{S^{\chi}_{M_1^{{\even}},M_2^{{\even}}}}\right)^2\\
&=\dfrac{1}{{\sin^2\left( \dfrac{2 \pi (k+2)}{2k+3} \right)}}
\sum_{1\leq r \leq 2k+2} 2\cdot{\sin^2\left( \dfrac{2 r \pi (k+2)}{2k+3} \right)}
= \frac{2k+3}{\sin^2\left(\frac{\pi}{2k+3} \right)}.
\end{align*}

\end{proof}	

\begin{cor}
The Frobenius-Perron dimension of the category of local $\eL$-modules, i.e., $\FP(\rep^0\eL)$
and the one of the sub tensor category given by the $M_r^{\text{even/odd}}$, i.e., $\FP(\cS_k)$  coincide.
\end{cor}
\begin{proof}
Immediate from Lemmas \ref{lem:sin2}, \ref{lem:FPrepA}, \ref{lem:FPCk} and  \ref{lem:FPMrevenodd} and equation \eqref{eqn:dimLkeven}.
\end{proof}
\begin{cor}
The $M_r^{\text{even/odd}}$ for $1\leq r\leq 2k+2$ form a complete list of simple local $\eL$-modules. 
\end{cor}

\begin{cor}
Let $1\leq r\leq 2k+2$. The $M_r$ for $r$ odd form a complete list (up to parity and isomorphisms) of simple local $\cL$-modules and for $r$ even of twisted ones. 
\end{cor}

\begin{rem}
\cite[Thm.\ 2.67]{CKM} says that the induction functor is a braided tensor functor when restricted to those objects that induce to local modules for an algebra object (if inducing to modules for a super algebra, one needs to start with an auxiliary super-category).
In our example of $\eF$ this is the full abelian  subcategory of $\cC$ formed by $\mathcal L_{1,0} \otimes V_{1,r}$ for $r=1, \dots, 2k+2$. Denote this category by $\mathcal C^{{\even}}_{k,0}$.
Further, let  $\mathcal{S}^{{\even}}$ be the full abelian  subcategory formed by $M_r^{{\even}}$-modules inside $\eL$-modules
and let $\mathcal{S}^{\text{super}}$ be the  full abelian subcategory of  $\mathcal{S}^{{\even}}$ formed by all $M_r^{{\even}}$ with $r$ odd. 

Then the same argument as in the proof of \cite[Thm.\ 5.1]{OS} shows that the induction functor $\eF$ restricted to $\mathcal C^{{\even}}_{k,0}$ is fully faithful and thus this subcategory is braided equivalent to $\mathcal{S}^{{\even}}$.

The argument of \cite[Thm.\ 5.1]{OS} also works for extensions by superalgebra objects (using then also \cite[Lem.\ 2.61]{CKM}) and thus one also gets a braided equivalence between 
an auxiliary super category $\mathcal{S}\mathcal{S}^{\text{super}}$ (see \cite[Defn.\ 2.11]{CKM})
and the category of local $\cL$-modules. 
As a consequence of these two equivalences of braided tensor categories, we especially have that the category of local $\cL$-modules is braided equivalent to the super-category auxiliary to the
full abelian subcategory of the vertex tensor category of the Virasoro algebra $\Vir(p, \Delta/2)$ formed by the $V_{1, r}$ with $r$ odd. 
\end{rem}

\subsection{Fusion rules}\label{sec:fusion}

Using that the induction functor is a $P(z)$-tensor functor we get the fusion rules of $\cL$ for free, namely
\begin{align}
M_r \boxtimes_{\cL} M_{r'} &\cong  
\bigoplus_{r'' =1 }^{2k+2} N_{r, r'}^{r''} M_{r''}
\end{align}
which we have established in \eqref{eqn:Mrfus}.
Each $M_r$ decomposes as $M_r=M_r^{{\even}} \oplus M_r^{{\odd}}$ as $\eL$-module. As explained in \cite{CKM} one also assigns a parity to each \svoa{} module and we have two choices, we can either give $M_r^{{\even}}$ even parity and $M_r^{{\odd}}$ odd one or the other way around. Let us denote the first choice by $M_r^+$ and the second one by $M_r^-$. 
An intertwining operator is then called even if it respects the chosen parities and odd otherwise. The super dimension ($\sdim$) of a fusion rule is accordingly the difference of the dimension of even and odd intertwining operators of the given type. From our discussion of $\eL$, we have the $\sdim$  of type \[
\sdim {\left( \substack{  M^{\epsilon''}_{r''} \\ M^{\epsilon}_r \ \  M^{\epsilon'}_{r'}  }\right)}  = 
N^{-\qquad\quad (r'', \epsilon'')}_{(r, \epsilon), (r', \epsilon')} = \epsilon \epsilon' \epsilon'' {N_{r, r'}}^{r''}.
\]
The fusion coefficients can also be computed using Verlinde's formula, see \cite[Sec.\ 1.5]{CKM}.

\subsection{Modular transformations}\label{sec:modular}

The results of this section follow by applying \cite[Sec.\ 4.2.1]{CKM}.
Let $1\leq r,r'\leq 2k+2$. The character and super character of $M_r^+$ are defined as 
\[
 \ch^{\pm}[M^+_{r}](\tau, v) =  \ch[M^{{\even}}_{r}](\tau, v)\pm  \ch[M^{{\odd}}_{r}](\tau, v).
\]
Define the numbers
\begin{equation}
s_{r, r'}:= (-1)^{r+r'}\sqrt{\frac{1}{2k+3}} \sin\left(\frac{\pi r r' (k+2)}{(2k+3)}\right).
\label{eqn:SLeven}
\end{equation}
Then, \cite[Prop.\ 2.89]{CKM}  and \cite[Thm.\ 4.5]{KO}  give
\begin{align}
S&{}^\chi_{M_r^{{\even}}, M_{r'}^{{\even}} } = 
\dfrac{1}{\mathcal{D}(\rep^0\eL)}S^{\hopflink}_{M_r^{{\even}}, M_{r'}^{{\even}} }
= \dfrac{\dim_{\cC}(\eL)}{\mathcal{D}(\cC)}S^{\hopflink}_{M_r^{{\even}}, M_{r'}^{{\even}} }
\nonumber \\
&= \dfrac{\dim_{\cC}(\eL)}{\mathcal{D}(\cC)}
S^{\hopflink}_{\mathcal{L}_{1,0}\otimes V_{1,r}, \mathcal{L}_{1,0}\otimes V_{1,r'}}
=\dim_{\cC}(\eL)
S^{\chi}_{\mathcal{L}_{1,0}\otimes V_{1,r}, \mathcal{L}_{1,0}\otimes V_{1,r'}}
=s_{r, r'}\label{eqn:SMeMe}.
\end{align}
Using \cite[Lem.\ 4.26]{CKM} we get 
\begin{align}
S^\chi_{M_r^{{\even}}, M_{r'}^{{\odd}} } = \begin{cases} s_{r, r'} \qquad &r \ {\even} \\ -s_{r, r'} \qquad &r \ {\odd} \end{cases} 
\qquad  \text{and} \qquad 
S^\chi_{M_r^{{\odd}}, M_{r'}^{{\odd}} } = \begin{cases} s_{r, r'} \qquad &r+r' \ {\even} \\ -s_{r, r'} \qquad &r+r' \ {\odd} \end{cases} 
\label{eqn:SMMother}.
\end{align}
We then have for $r$ odd
\begin{equation}
\begin{split}
  \ch^+[M^+_r]\left(-\frac{1}{\tau}, \frac{z}{\tau}\right) &= e^{2\pi i k \frac{z^2}{\tau}} \sum_{r' \ {\even}} 2s_{r, r'}  \ch^-[M^+_{r'}](\tau, z) \\
\ch^-[M^+_r]\left(-\frac{1}{\tau}, \frac{z}{\tau}\right) &= e^{2\pi i k \frac{z^2}{\tau}}  \sum_{r' \ {\odd}} 2s_{r, r'}   \ch^-[M^+_{r'}](\tau, z)
\end{split}
\end{equation}
and for $r$ even
\begin{equation}
\begin{split}
 \ch^+[M^+_r]\left(-\frac{1}{\tau}, \frac{z}{\tau}\right) &=e^{2\pi i k \frac{z^2}{\tau}}  \sum_{r' \ {\even}} 2s_{r, r'}  \ch^+[M^+_{r'}](\tau, z) \\
\ch^-[M^+_r]\left(-\frac{1}{\tau}, \frac{z}{\tau}\right) &= e^{2\pi i k \frac{z^2}{\tau}}  \sum_{r' \ {\odd}} 2s_{r, r'}   \ch^+[M^+_{r'}](\tau, z).
\end{split}
\end{equation}
Assemble the $S^\chi$ into an $(4k+4)\times (4k+4)$-matrix with first $(2k+2)$ rows and columns corresponding to 
$M_r^{{\even}}$ and the remaining rows and columns with corresponding $M_r^{{\odd}}$.
By changing the basis of the Grothendieck ring so that $M_r^{{\even}}, M_r^{{\odd}}$
are replaced by $M_r^{{\even}}+ M_r^{{\odd}}$ and $M_r^{{\even}}- M_r^{{\odd}}$
as in \cite[Eqns.\ (4.9), (4.10)]{CKM} form the $(4k+4)\times (4k+4)$ matrices $\widetilde{S}$ and $\widetilde{S}^{-1}$. Following the notation
in \cite{CKM}, Verlinde's formula \cite[Eqn.\ (4.11)]{CKM} reads
\begin{equation}
\begin{split}
N^{+\qquad\quad (r'', +)}_{(r, +), (r', +)} &=  
\delta(r+r'+r'' =1 \mod 2)\sum\limits_{t \ {\even}} 
\dfrac{\widetilde{S}_{r, t}\cdot {\widetilde{S}}_{r' ,t}\cdot (\widetilde{S}^{-1})_{t, r''}}{\widetilde{S}_{1, t}}\\
N^{-\qquad\quad (r'', +)}_{(r, +), (r', +) }&=\delta(r+r'+r'' =1 \mod 2)
\sum\limits_{t \ {\odd}} \dfrac{\widetilde{S}_{r, t}\cdot\widetilde{S}_{r' ,t}\cdot (\widetilde{S}^{-1})_{t, r''}}{\widetilde{S}_{1, t}}.	
\end{split}
\end{equation}

\section{Parafermions}\label{sec:para}

$\cL$ contains the lattice \voa{} $V_L$ of the lattice $L=\sqrt{2k}\mathbb Z$ as vertex operator subalgebra.
Let $C_k := \text{Com}\left( V_{L}, \cL\right)$ the parafermion coset. It is rational by 
\cite[Cor.\ 4.13]{CKLR} and since
$C_k = \text{Com}\left( V_{L}, \eL\right)$.
Every even weight module of $L_k(\mathfrak{sl}_2)$ is graded by $2L'$ and every odd weight module is graded by 
$L'\backslash2L'$, it follows especially that $\cL$ is graded by $L'$ and hence, for $1\leq r \leq 2k+2$,
\[
M_r \cong \bigoplus_{\lambda \in L'/L} V_{\lambda+L} \otimes C_{\lambda, r}. 
\]
We restrict to local modules $M_r$, that is $r$ is odd. 
The results of \cite{CKLR, CKL, CKM} together, (see \cite[Sec.\ 4.3.1]{CKM}), now imply that all $C_{\lambda, r}$ are simple $C_k$-modules. The grading lattice being the full dual lattice $L'$ implies (\cite[Thm.\ 4.39]{CKM}) that $C_{\lambda, r} \cong C_{\nu, s}$ if and only if $\lambda=\nu$ and $r=s$. 
Moreover, these are all inequivalent simple modules, this follows from Thm.\ 4.3, Lem.\ 4.9 and the proof of Thm.\ 4.12 of \cite{CKLR}. 
We remark that these Theorems are formulated for \voas{} and not \svoas. Modifying the proofs for the \svoa{} setting is not difficult. 
The easier route is however to apply the Theorems of \cite{CKLR, CKL, CKM} to $\eL$ in which case each simple $C_k$-module appears as a submodule of two distinct $\eL$-modules. One of these two $\eL$-modules then lifts to a local $\cL$-module and the other one to a twisted one. 

Using the induction functor (\cite[Thm.\ 4.41]{CKM}) one deduces that 
\[
C_{\nu, r} \boxtimes C_{\lambda, r'} \cong \bigoplus_{r'' =1 }^{2k+2} {N_{r, r'}}^{r''} C_{\lambda+\nu, r''}.
\]
From \cite[Thm.\ 4.42]{CKM}, we have that:
\begin{align*}
\ch^+[M_r^+](\tau, u) &= 
\sum_{\nu\in L'/L}\dfrac{\theta_{L+\nu}(u,\tau)}{\eta(\tau)}
\ch[C_{\nu,r}](\tau),\\
\ch^-[M_r^+](\tau, u) &= 
\sum_{\nu\in L'/L}\dfrac{\theta_{L+\nu}(u,\tau)}{\eta(\tau)} 
(-1)^{\delta(\,\nu \not\in (2L')/L\,)}\ch[C_{\nu,r}](\tau)
,\\
\ch[C_{\nu,r}](\tau)&=\dfrac{\eta(\tau)}{2k\cdot \theta_{L+\nu}(0,\tau)}
\sum\limits_{\gamma\in L'/L }e^{-\pi i \langle \nu,\gamma \rangle}
\ch^+[M_r^+](\tau,\gamma),\\
&=(-1)^{\delta(\,\nu \not\in (2L')/L\,)}\dfrac{\eta(\tau)}{2k\cdot \theta_{L+\nu}(0,\tau)}
\sum\limits_{\gamma\in L'/L }e^{-\pi i \langle \nu,\gamma \rangle}
\ch^-[M_r^+](\tau,\gamma),
\end{align*}
for the $T$ transformations, we get:
\begin{align*}
e^{\pi i \left(\langle \lambda,\lambda\rangle-\frac{1}{12}\right)} T_{C_{\lambda,r}} &= T_{M_r},
\end{align*}
and finally, we can calculate the $S$-transformations by inducing up to $\cL^{{\even}}$-modules, 
and then using equations \eqref{eqn:SLeven}, \eqref{eqn:SMeMe}, \eqref{eqn:SMMother} as in \cite[Thm.\ 4.42]{CKM}:
\begin{align*}
S^{\chi, C}_{C_{\lambda,r_1}, C_{\mu,r_2}}\cdot e^{2\pi i \langle \lambda, \mu \rangle}\cdot \sqrt{k/2} &= 
\begin{cases}
s_{r,r'} & \text{ if } \lambda,\mu \in 2L'/L,\\
(-1)^rs_{r,r'} & \text{ if } \lambda\in 2L'/L, \mu\not\in 2L'/L \text{ or } \mu\in 2L'/L, \lambda\not\in 2L'/L, \\
(-1)^{r+r'}s_{r,r'} & \text{ if } \lambda,\mu\not\in 2L'/L.
\end{cases}
\end{align*}

\bibliographystyle{alpha}

\end{document}